\documentclass[francais]{smfart}
\usepackage{epsfig,amsmath,a4,amssymb,amsfonts,amsbsy,xy,latexsym,epsf,pstricks,multicol,accents}
\xyoption{all}
\usepackage[french]{babel}
\theoremstyle{plain}

\let\set\mathbb
\def\<<{\leavevmode
  \raise0.28ex\hbox{$\scriptscriptstyle\langle\!\langle$}\nobreak
  \hskip -.6pt plus.3pt minus.2pt\,}
\def\>>{\,\nobreak\hskip -.6pt plus.3pt minus.2pt
  \raise0.28ex\hbox{$\scriptscriptstyle\rangle\!\rangle$}}

\newcommand{\VEC}[1]{\overrightarrow {#1}}

\def\pf{{\oplus _\cF}}

\def\U{{\set U}}
\def\Zp{{\ZZ_p}}
\def\Zps{{\ZZ_p^*}}
\def\Qp{{\QQ_p}}
\def\Qps{{\QQ_p^*}}
\def\V{{\set V}}
\def\Log{\mathop{\rm{Log}}\nolimits }
\def\Logf{\mathop{\rm{Log}_\cF}\nolimits }
\def\expf{\mathop{\rm{exp}_\cF}\nolimits }

\def\SL2Z{{{\rm SL}_2(\ZZ)}}

\def\clk{{{\rm  Pic}(\OO)}}

\def\bG{{\rm G}}

\def\AA{{\set A}}
\def\Aut{\mathop{\rm{Aut}}\nolimits }

\def\FF{{\set F}}

\def\Fp{{\FF _p}}

\def\Fq{{\FF _q}}
\def\Fqs{{\FF^*_q}}

\def\Gal{\mathop{\rm{Gal}}\nolimits }

\def\Ker{\mathop{\rm{Ker}}\nolimits }
\def\NN{{\set N}}
\def\PP{{\set P}}
\def\PP{{\set P}}

\def\Pic{\mathop{\rm{Pic}}\nolimits }
\def\QQ{{\set Q}}

\def\TT{{\set T}}

\def\ZZ{{\set Z}}
\def\agot{{\mathfrak a}}
\def\fgot{{\mathfrak f}}

\def\Gm{\mathop{\rm G_m}\nolimits }
\def\Ga{\mathop{\rm G_a}\nolimits }
\def\GL{\mathop{\rm GL}\nolimits }
\def\GG{\mathop{\rm Gal}\nolimits }
\def\End{\mathop{\rm End}\nolimits }

\def\OO{{\mathcal O}}

\def\cO{{\mathcal O}}
\def\cB{{\mathcal B}}
\def\cC{{\mathcal C}}

\def\cE{{\mathcal E}}
\def\cF{{\mathcal F}}
\def\cG{{\mathcal G}}

\def\cN{{\mathcal N}}
\def\cO{{\mathcal O}}

\def\cT{{\mathcal T}}
\def\cU{{\mathcal U}}

\def\cX{{\mathcal X}}

\def\lgot{{\mathfrak l}}

\def\sc{\scriptscriptstyle }

\newtheorem{theoreme}{Th{\'e}or{\`e}me}

\newtheorem{problem}{Probl{\`e}me}

\author[Jean-Marc Couveignes]{Jean-Marc  Couveignes\thanks{L'auteur est
    soutenu par le fond national pour la science (ACI NIM), et par le Centre
    d'\'electronique de l'Armement (CELAR, DGA)}}
\address{Institut de math{\'e}matiques de Toulouse, Universit{\'e} de Toulouse et CNRS}
\email{couveig@univ-tlse2.fr}
\urladdr{http://www.picard.ups-tlse.fr/~couveig}
\title{Quelques  math{\'e}matiques de la cryptologie {\`a} cl{\'e}s publiques}
\alttitle{A few   mathematical tools for public key cryptology}

\begin{document}

\frontmatter

\begin{abstract}
Cette note pr{\'e}sente quelques d{\'e}veloppements math{\'e}matiques plus ou
moins r{\'e}cents de la
cryptologie {\`a} cl{\'e}s publiques. 
\end{abstract}

\begin{altabstract}
I present  examples of mathematical objects that are of interest for 
public key cryptography.
\end{altabstract}


\altkeywords{jacobian variety,  algorithmic complexity, public key cryptology,
commutative algebraic group, algebraic curve, finite field}

\bibliographystyle{smfplain}
\maketitle
\let\languagename\relax
\mainmatter

\tableofcontents

La cryptologie {\`a} cl{\'e}s publiques a mobilis{\'e} depuis son
invention  des math{\'e}matiques plus ou moins {\'e}l{\'e}mentaires : arithm{\'e}tique des
congruences, th{\'e}orie alg{\'e}brique des nombres, g{\'e}om{\'e}trie et cohomologie des groupes
alg{\'e}briques, th{\'e}orie des graphes, probabilit{\'e}s discr{\`e}tes, complexit{\'e}
algorithmique, etc. 
Il n'est pas
toujours facile de discerner une ligne directrice dans ces d{\'e}veloppements et
c'est
une des difficult{\'e}s  du domaine. Je me bornerai donc {\`a} pr{\'e}senter quelques
id{\'e}es et situations typiques, sans aucune pr{\'e}tention {\`a} l'exhaustivit{\'e}.

La section \ref{section:H2S} d{\'e}crit deux protocoles classiques de la cryptographie {\`a}
cl{\'e}s publiques dans le  cadre g{\'e}n{\'e}ral de l'action d'un groupe sur un ensemble
fini. Le logarithme discret offre une exemple de cette situation. Les groupes
utilis{\'e}s sont le plus souvent des groupes de points rationnels
d'un groupe alg{\'e}brique commutatif sur un corps fini. Je montre dans la
section  \ref{section:LOG} quel genre de propri{\'e}t{\'e}s on attend (ou on redoute) 
d'un groupe alg{\'e}brique dans ce contexte. La section \ref{section:isog} met en jeu
non plus un groupe alg{\'e}brique mais une cat{\'e}gorie de groupes,  leurs
morphismes, et les graphes  qui s'en d{\'e}duisent.

\section{Espaces homog{\`e}nes difficiles}\label{section:H2S}

Dans cette section nous  d{\'e}crivons 
une famille de probl{\`e}mes calculatoires que nous appelons
espaces homog{\`e}nes difficiles (EHD). Nous montrons
que cette notion offre un cadre  naturel
{\`a} nombre de protocoles  fondamentaux
de la cryptologie {\`a} cl{\'e} publique, pour le chiffrement,
l'identification, et l'{\'e}change de cl{\'e}  par exemple.

Le probl{\`e}me du logarithme discret (dans un  groupe
multiplicatif ou une courbe elliptique sur un corps fini) 
fournit un exemple d'EHD. Mais il existe
bien d'autres EHD. Ceux que nous pr{\'e}sentons dans la section
\ref{section:isog}
proviennent
de la multiplication complexe des courbes elliptiques.

Dans le paragraphe \ref{subsection:H2Sdef} 
nous d{\'e}finissons les espaces
homog{\`e}nes difficiles. Nous expliquons au paragraphe
\ref{subsection:DL} pourquoi le logarithme discret est un cas
particulier d'espace homog{\`e}ne.
Nous pr{\'e}sentons dans le paragraphe \ref{subsection:DHM}
un protocole d'{\'e}change de cl{\'e}s de type  Diffie-Hellman-Merkle
dans le contexte des EHD. 
Nous d{\'e}crivons  de m{\^e}me dans le paragraphe \ref{subsection:schnorr}
le  protocole de Schnorr pour la preuve de connaissance
sans apport d'information, 
dans le  cadre g{\'e}n{\'e}ral
et naturel  des EHD.

\subsection{D\'efinition d'un espace homog{\`e}ne difficile}\label{subsection:H2Sdef}

Soit $G$ un groupe fini commutatif.
Un espace homog{\`e}ne $H$ pour $G$ est un ensemble fini non vide 
$H$ muni d'une action transitive et libre de $G$.
Donc le cardinal de $H$ est {\'e}gal {\`a} celui de $G$. On note
 $S=\# H=\# G$. On appelle {\it points} les {\'e}l{\'e}ments de $H$
et {\it vecteurs} les {\'e}l{\'e}ments de $G$. Un exemple
naturel : $H$ est un espace affine et $G$ l'espace vectoriel sous-jacent.

Si   $h_1$ et $h_2$  sont deux points,  il existe
un  unique vecteur $g$ tel que $g.h_1=h_2$. On note $\VEC{h_1h_2}$
ce vecteur. 

\'Etant donn{\'e} un espace homog{\`e}ne, on consid{\`e}re
une s{\'e}rie de probl{\`e}mes calculatoires.

On suppose que les {\'e}l{\'e}ments de $G$ et de $H$ sont repr{\'e}sent{\'e}s par
des cha{\^\i}nes de caract{\`e}res de longueur polynomiale en $\log S$.

On doit {\^e}tre capable de calculer efficacement
la loi de composition et
l'inversion dans le groupe $G$ et de tester l'{\'e}galit{\'e} 
de deux {\'e}l{\'e}ments de ce groupe. Autrement dit, on veut que
le groupe $G$ soit calculatoire.

\begin{problem}[Op\'erations dans le groupe $G$]\label{pb:opG}
\'Etant  donn{\'e}s deux vecteurs $g_1$ et $g_2$,
d{\'e}cider s'ils sont {\'e}gaux, calculer l'inverse  $g_1^{-1}$
de $g_1$ et le produit 
$g_1g_2$.
\end{problem}

Il faut aussi pouvoir choisir des {\'e}l{\'e}ments
al{\'e}atoires dans $G$.

\begin{problem}[Vecteur al{\'e}atoire]\label{pb:alea}
Choisir un vecteur $g$ dans $G$ avec une distribution
(presque) uniforme.
\end{problem}

On souhaite aussi r{\'e}soudre efficacement les probl{\`e}mes {\'e}l{\'e}mentaires
suivants concernant l'action de $G$ sur $H$ :

\begin{problem}[Action de $G$ sur $H$]\label{pb:opH}
\'Etant donn{\'e}s deux points  $h_1, h_2\in H$  et un
vecteur $g\in G$, d{\'e}cider si $h_1=h_2$, et calculer $g.h_1$.
\end{problem}

Notons que si l'on applique un vecteur al{\'e}atoire
(avec distribution uniforme)
$g$ {\`a} un point fixe $h_0$, on obtient un
point al{\'e}atoire (avec distribution uniforme).

On dit que l'espace homog{\`e}ne 
est calculatoire si l'on dispose
d'un algorithme probabiliste 
polynomial en $\log S$ pour
r{\'e}soudre les  probl{\`e}mes \ref{pb:opG}, \ref{pb:alea} et \ref{pb:opH}.
Cela sous-entend que l'on consid{\`e}re, non pas un espace homog{\`e}ne,  mais une 
famille infinie d'espaces homog{\`e}nes.

Venons en maintenant {\`a} des propri{\'e}t{\'e}s plus subtiles.

Souvenons nous qu'il existe un unique vecteur $\VEC{h_1h_2}$ 
qui envoie $h_1$ sur $h_2$ :
$$\VEC{h_1h_2}.h_1=h_2.$$

On peut souhaiter calculer ce vecteur.

\begin{problem}[Diff{\'e}rence de deux points]\label{pb:vec}
\'Etant donn{\'e}s  $h_1, h_2\in H$ trouver  $g\in G$
tel que  $g.h_1=h_2$.
\end{problem}

Un probl{\`e}me de m{\^e}me nature 
est de compl{\'e}ter un parall{\'e}logramme.

\begin{problem}[Compl{\'e}tion d'un  parall{\'e}logramme]\label{pb:CDH}
\'Etant donn{\'e}s trois  points   $h_1, h_2, h_3 \in H$, calculer 
l'unique point  $h_4$ tel que  $\VEC{h_1h_2}=\VEC{h_3h_4}$.
\end{problem}

Ce  $h_4$ n'est autre que  $\VEC{h_1h_2}.h_3$. Donc le probl{\`e}me \ref{pb:CDH}
est plus facile que le probl{\`e}me \ref{pb:vec}.

\vskip 0.3cm

On s'int{\'e}resse aux espaces homog{\`e}nes calculatoires 
pour lesquels  les probl{\`e}mes \ref{pb:vec} et \ref{pb:CDH}
sont difficiles. Cela signifie qu'il n'existe pas de machine
de Turing probabiliste qui r{\'e}solve l'un ou l'autre
de ces probl{\`e}mes  en temps polynomial en $\log S$.

De tels espaces homog{\`e}nes sont appel{\'e}s espaces homog{\`e}nes difficiles
(EHD).

On pourrait consid{\'e}rer un autre probl{\`e}me 

\begin{problem}[V{\'e}rification d'un parall{\'e}logramme]\label{DDH}
\'Etant donn{\'e}s quatre points  $h_1,h_2,h_3,h_4$ dans  $H$, dire si
$\VEC{h_1h_2}=\VEC{h_3h_4}$.
\end{problem}

Si ce dernier probl{\`e}me est difficile
on dit que l'espace homog{\`e}ne est 
tr{\`e}s difficile (EHTD).

Supposons que  $G=k^d$ est un espace vectoriel sur un corps
fini $k$ et $H=\AA^d(k)$ l'espace affine associ{\'e}. C'est un espace
calculatoire.
Les vecteurs et les points sont d{\'e}crits pas leurs coordonn{\'e}es et  l'action de 
$g=(x_1,\ldots, x_d)$  sur $h=(a_1,\ldots, a_d)$ se calcule
au prix de $d$ additions dans $k$.

Ce n'est pas un espace homog{\`e}ne difficile car si $h=(a_1,\ldots, a_d)$ 
et $k=(b_1,\ldots, b_d)$ alors $\VEC{hk}=(b_1-a_1,\ldots, b_d-a_d)$ 
se calcule au prix de $d$ soustractions dans $k$.

\subsection{Le logarithme discret}\label{subsection:DL}
Un premier candidat EHD int{\'e}ressant  est fourni par le
probl{\`e}me du logarithme discret.

Soit   $C$ un groupe cyclique d'ordre 
$n$ et soit $c$ un g{\'e}n{\'e}rateur de $C$.

Notons $G$ le groupe des automorphismes de $C$. Un {\'e}l{\'e}ment
$g$ de $G$ envoie $c$ sur $g(c)=c^a$ o{\`u} $a$ est un entier
premier {\`a} $n$. L'application $g\mapsto a$ est un isomorphisme
de $G$ sur $(\ZZ/n\ZZ)^*$.

Soit  $H$ l'ensemble des g{\'e}n{\'e}rateurs de $C$. Alors
$\# H= \# G = \phi (n)$ et $G$ agit simplement transitivement
sur $H$.

On suppose que $C$ est un groupe calculatoire, que son ordre $n$
est connu, et que la factorisation
de $n$ en produit de facteurs premiers est connue elle aussi.
Alors
on dispose d'algorithmes polynomiaux en $S=\phi(n)$ pour
calculer dans $G=\ZZ/\phi(n)\ZZ$. R{\'e}soudre le probl{\`e}me
\ref{pb:opH} revient {\`a}  calculer $c^a$ pour $c$ un g{\'e}n{\'e}rateur
de $C$ et $a$ un entier entre $1$ et $n$.

Un algorithme na\"\i f calculerait successivement $c$, $c^2$, $c^3$,
$c^4$, \ldots, $c^{a-1}$, $c^a$,  ce qui requiert $a-1$ op{\'e}rations.
Ce n'est pas satisfaisant car on souhaite calculer $c^a$ en temps
$\log S$.

L'algorithme utilis{\'e} est connu sous le nom d'exponentiation rapide.
On calcule $c_0=c$, $c_1=c_0^2=c^2$, $c_2=c_1^2=c^4$, $c_3=(c_2)^2=c^8$,  \ldots, $c_x=
c^{2^x}$
o{\`u} $2^x$ est la plus grande puissance de $2$ inf{\'e}rieure ou {\'e}gale 
{\`a} $a$. On {\'e}crit alors l'exposant $a$ en base $2$ 
soit $a=\sum_{1\le k\le
  x}\epsilon_k2^k$ et on v{\'e}rifie que $c^a=\prod_{1\le k\le
  x}c_k^{\epsilon_k}$. Au total le calcul de  $c^a$ n'a pas requis
plus de $2\log_2 a$ op{\'e}rations dans $C$.

En revanche, si $C$ est  un  groupe quelconque, alors
on ne sait pas, en g{\'e}n{\'e}ral,  r{\'e}soudre efficacement les 
probl{\`e}mes $\ref{pb:vec}$
et $\ref{pb:CDH}$. Par exemple, le probl{\`e}me
$\ref{pb:vec}$ dans ce contexte est le suivant : {\'e}tant donn{\'e}s
deux g{\'e}n{\'e}rateurs $c$ et $d$ de $C$, trouver un entier $k$
tel que $d=c^k$. Cet entier $k\in \ZZ/n\ZZ$ est appel{\'e} logarithme
discret de $d$ en base $c$ et not{\'e} $\log_c(d)$.

On observe que $\log_c$ s'{\'e}tend en une application $\log_c : C\rightarrow
\ZZ/n\ZZ$ qui est un isomorphisme
de groupe. C'est l'application r{\'e}ciproque de l'exponentiation 
de base $c$ not{\'e}e $exp_c : \ZZ/n\ZZ\rightarrow C$ et d{\'e}finie
par $\exp_c(k)=c^k$.

Le probl{\`e}me du {\it logarithme discret}  a un sens
pour tout groupe cyclique et c'est  un cas
particulier d'espace homog{\`e}ne.

Nous avons  l{\`a} un  premier exemple  
d'espace homog{\`e}ne difficile. 
Il n'existe pas, en effet,  d'algorithme g{\'e}n{\'e}rique pour calculer 
les logarithmes discrets\footnote{Un algorithme g{\'e}n{\'e}rique est 
un algorithme qui
  n'utilise pas d'autres propri{\'e}t{\'e}s de $C$ que l'existence d'une loi de
  groupe. Shoup a montr{\'e} qu'un tel algorithme ne peut pas calculer le
  logarithme discret en temps $o(\sqrt P)$ o{\`u} $P$ est le plus grand facteur
premier de l'ordre $n$ de $C$.}. 

Une petite difficult{\'e} subsiste : il n'existe pas de groupe
g{\'e}n{\'e}rique. Rien n'interdit {\`a} un algorithme d'utiliser des propri{\'e}t{\'e}s
particuli{\`e}res au groupe utilis{\'e} en pratique. 
On ne conna{\^\i}t pas d'algorithme
pour r{\'e}soudre le logarithme discret en temps polynomial en $\log n$
dans les groupes multiplicatifs 
de corps finis et on suppose qu'il n'en existe pas. Cela ne prouve
pas cependant 
qu'il n'en existe pas. Il est g{\'e}n{\'e}ralement admis cependant que
les espaces homog{\`e}nes correspondants
sont difficiles et m{\^e}me tr{\`e}s difficiles.

Ce premier exemple d'EHD est aussi un exemple (tout aussi hypoth{\'e}tique)
de fonction asym{\'e}trique. Cel{a} signifie que les deux fonctions
$\exp_c$ et $\log_c$ sont 
r{\'e}ciproques l'une de l'autre, que $\exp_c$ se calcule en temps polynomial
(gr{\^a}ce {\`a} l'algorithme d'exponentiation rapide) mais que
$\log_c$ ne se calcule pas en temps polynomial.

On trouve dans \cite{Salomaa} une introduction {\`a} la cryptologie
asym{\'e}trique. Le livre \cite{Goldreich} est une introduction g{\'e}n{\'e}rale
aux concepts de la cryptologie moderne. Le
livre   \cite{Rolland} est un  trait{\'e} g{\'e}n{\'e}ral
et r{\'e}cent  de 
cryptologie.

\subsection{\'Echange de cl{\'e}}\label{subsection:DHM}

Nous pr{\'e}sentons dans ce paragraphe un exemple de protocole cryptographique
qui repose sur un EHD. Il s'agit de la transposition {\'e}vidente  du protocole 
de  Diffie-Hellman-Merkle dans le contexte des EHD.

On suppose qu'Alice et Bob communiquent par un canal
non s{\'e}curis{\'e} (tout le monde 
 peut entendre ou lire l'int{\'e}gralit{\'e}
de leurs messages). Au d{\'e}but du protocole ils ne partagent
aucun secret. \`A l'issue du protocole, ils ont un secret
commun, c'est-{\`a}-dire une information connue d'eux seuls. Cette information
pourra leur servir de cl{\'e} secr{\`e}te pour des {\'e}changes ult{\'e}rieurs.

Voici comment ils proc{\`e}dent. Ils conviennent publiquement
d'un EHD $H$. Dans tout ce qui suit, le mot al{\'e}atoire sous-entend 
que la distribution est uniforme (ou tr{\`e}s proche de la distribution
uniforme).

\begin{enumerate}
\item  Alice choisit un point  al{\'e}atoire 
 $h_0$ dans  $H$ et un vecteur al{\'e}atoire 
 $g_1$ dans  $G$. Elle applique $g_1$ {\`a} $h_0$
et calcule   $h_1=g_1.h_0$. Elle envoie le couple
de points  $(h_0,h_1)$ {\`a} Bob.

\item Bob choisit un vecteur al{\'e}atoire 
 $g_2$ dans  $G$ et l'applique {\`a} $h_0$. Il calcule donc
 $h_2=g_2.h_0$. Il envoie  $h_2$ {\`a} Alice. La 
cl{\'e} secr{\`e}te est  $K=g_2.h_1$.

\item Alice calcule la cl{\'e} secr{\`e}te {\`a} partir des informations
dont elle dispose : elle applique $g_1$ {\`a} $h_2$. En effet   $K=g_2.h_1=g_1.h_2$.

\end{enumerate}

La commutativit{\'e} de $G$ joue ici un r{\^o}le essentiel. Alice et Bob conviennent
d'une origine $h_0$ publique. Ils choisissent chacun un vecteur
et construisent ensemble un parall{\'e}logramme. Trois sommets  $h_0$,
$h_1$, $h_2$ du 
parall{\'e}logramme sont publics mais le quatri{\`e}me $K$
 est connu d'eux seuls. Alice conna{\^\i}t $K$ car elle a choisi
le cot{\'e} $g_1$ et elle a re{\c c}u de Bob  le sommet $h_2$. 
Bob conna{\^\i}t $K$ car il a choisi
le cot{\'e} $g_2$ et il a re{\c c}u d'Alice  le sommet $h_1$. 
Un observateur {\'e}tranger {\`a} cet {\'e}change voit les trois sommets $h_0$,
$h_1$ et $h_2$ mais aucun des cot{\'e}s du parall{\'e}logramme. Il doit
donc compl{\'e}ter le parall{\'e}logramme pour violer le secret commun
{\`a} Alice et Bob. Et on suppose que ce calcul est trop difficile (il n'existe
pas d'algorithme polynomial en temps pour le r{\'e}soudre).

\subsection{Preuve de connaissance sans apport d'information}\label{subsection:schnorr}

Nous pr{\'e}sentons dans ce paragraphe  le protocole de Schnorr dans le cadre des 
EHD.

On suppose qu'Alice a choisi un EHD $H$ et un point 
$h_0$,  ainsi qu'un vecteur al{\'e}atoire $g_A$. Elle applique
le vecteur $g_A$ au point  $h_0$ et obtient le point $h_A=g_A.h_0$.

Elle publie $G$, $H$, $h_0$ et $h_A$ et garde $g_A$ secret. Le secret d'Alice
est donc le vecteur $g_A=\VEC{h_0h_A}$.

Alice veut prouver {\`a} Bob qu'elle conna{\^\i}t ce vecteur,  sans le divulguer.

\begin{enumerate}

\item  Alice choisit un vecteur  $g_r\in G$ al{\'e}atoire et calcule
$g_r.h_A=h_r$. Elle envoie  $h_r$ {\`a} Bob.

\item  Bob tire {\`a} pile ou face et envoie le r{\'e}sultat $\epsilon \in \{0,1\}$
  {\`a} Alice.

\item  Si  $\epsilon =0$ alors  Alice envoie  $g_p=g_r$ {\`a}
Bob. Sinon elle  envoie $g_p=g_rg_A$. 

\item  Bob v{\'e}rifie que  $g_p.h_A$ est {\'e}gal {\`a}  $h_r$ (si 
$\epsilon =0$) ou $g_p.h_0=h_r$ (si $\epsilon =1$).

\end{enumerate}

Le protocole construit un triangle $h_0$, $h_A$, $h_r$. Pour prouver
qu'elle conna{\^\i}t le cot{\'e} $\VEC{h_0h_A}$,  Alice prouve qu'elle conna{\^\i}t les deux
autres cot{\'e}s du triangle. Selon la valeur de $\epsilon$, Bob lui demandera
de d{\'e}voiler $\VEC{h_0h_r}$  ou $\VEC{h_Ah_r}$. Comme elle ne sait pas laquelle
de ces deux questions lui sera pos{\'e}e, elle doit conna{\^\i}tre la r{\'e}ponse au deux
questions.

Si elle ne conna{\^\i}t pas $\VEC{h_0h_A}$ elle ne peut conna{\^\i}tre {\`a} la fois
$\VEC{h_0h_r}$  et  $\VEC{h_Ah_r}$. Donc elle est prise en d{\'e}faut
par la question de Bob avec  une probabilit{\'e}  $\ge \frac{1}{2}$.

On  r\'ep{\`e}te le protocole un nombre suffisant de fois pour que Bob se convainque
qu'Alice
 conna{\^\i}t bien $\VEC{h_0h_A}$.

De son cot{\'e}, Bob n'apprend aucune information sur le secret d'Alice. Car tous 
les
vecteurs qui lui sont communiqu{\'e}s ont une extr{\'e}mit{\'e} al{\'e}atoire. Un observateur
ext{\'e}rieur
au protocole n'apprend rien lui non plus. Alice n'a donc pas besoin de r{\'e}v{\'e}ler
quoi que ce soit de son secret pour prouver {\`a} Bob qu'elle le conna{\^\i}t. On
mesure
l'avantage de cette m{\'e}thode sur le classique {\'e}change de mots de passe.

Le protocole d'{\'e}change de cl{\'e}s  a {\'e}t{\'e} publi{\'e} par Diffie et Hellman en 1976.
Le travail de Merkle sur cette question a jou{\'e} un r{\^o}le important dans
ce domaine. Il semble bien que ce protocole  ait {\'e}t{\'e} d{\'e}couvert ant{\'e}rieurement
par Williamson dans le cadre de son travail pour les services secrets
britanniques (raison pour laquelle il ne l'a pas publi{\'e}).

\section{Logarithmes discrets et groupes alg{\'e}briques}\label{section:LOG}

On cherche des groupes finis, cycliques, o{\`u} le logarithme discret soit
difficile.
Le groupe additif d'un corps fini n'est pas un bon candidat. Supposons par
exemple
que $C=(\ZZ/p\ZZ,+)$ pour $p$ premier et soient $g$ et $h$ deux r{\'e}sidus modulo
$p$. On suppose que $g$ engendre $C$. Donc $g$ est non nul modulo $p$.
Le logarithme discret de $h$ en base $g$ est l'entier $k$ tel que $h=kg$. Donc 
$k=\frac{h}{g}\bmod p$ se calcule {\`a} l'aide de l'algorithme d'Euclide en temps
$\le (\log p)^{2+o(1)}$ et m{\^e}me plus vite si l'on  a recours {\`a} des
algorithmes rapides. 

Le groupe le plus souvent utilis{\'e} est le groupe multiplicatif $\Gm (\Fq)=\Fqs$
 d'un corps fini $\Fq$. Les algorithmes connus  les plus rapides pour
 calculer les logarithmes discrets dans de tels groupes ont une complexit{\'e} 
de $\exp(\log(q)^{\frac{1}{3}+o(1)})$. Ce ne sont donc pas des algorithmes
polynomiaux. La contribution de Joux et Lercier {\`a} ce volume
est enti{\`e}rement consacr{\'e}e {\`a} cette question.

D'autres groupes alg{\'e}briques 
sont utilis{\'e}s depuis peu. Il s'agit principalement des courbes
elliptiques. Mais on a aussi sugg{\'e}r{\'e} 
l'utilisation de jacobiennes de courbes de genre sup{\'e}rieur (surtout le genre
deux). Les tores alg{\'e}briques sont  aussi l'objet d'{\'e}tudes
approfondies. Nous les pr{\'e}sentons  dans le paragraphe \ref{subsection:tores}.

Les meilleurs algorithmes connus pour calculer le logarithme discret dans le
groupe
des points d'une courbe elliptique  sur le corps {\`a} $q$
{\'e}l{\'e}ments, sont des algorithmes g{\'e}n{\'e}riques et ont donc une
complexit{\'e} 
en $\Omega(\sqrt P)$ o{\`u} $P$ est le plus grand facteur premier de l'ordre de la
courbe.  On en d{\'e}duit g{\'e}n{\'e}ralement que les cryptosyst{\`e}mes {\`a} bases
de courbes elliptiques peuvent atteindre un m{\^e}me niveau de s{\'e}curit{\'e} que ceux
bas{\'e}s sur les groupes multiplicatifs, avec une taille de cl{\'e} (la taille de
$q$) beaucoup plus petite (disons $200$ bits au lieu de $2000$ bits).

Il convient de souligner qu'il existe des instances faibles du logarithme
discret, tant pour les corps finis que pour les courbes elliptiques. Plus
pr{\'e}cis{\'e}ment,
il existe des familles infinies de groupes multiplicatifs et de courbes
elliptiques pour lesquels on dispose d'algorithmes polynomiaux de calcul du
logarithme discret. C'est {\'e}vident pour les groupes $C$ dont l'ordre $n$ n'a
pas de grand facteur premier. On donne des  exemples moins triviaux  au paragraphe
\ref{subsection:araki}.

\subsection{Un aper{\c c}u de l'utilisation des tores en
  cryptographie}\label{subsection:tores}

\subsubsection{Rappels sur les groupes alg{\'e}briques commutatifs}

Tout groupe alg{\'e}brique affine connexe de dimension
$1$ sur un corps alg{\'e}briquement clos est soit le groupe additif $\Ga$, soit le groupe  multiplicatif
$\Gm $.

Il y a deux familles importantes de groupes alg{\'e}briques.
Une  {\it  vari{\'e}t{\'e} ab{\'e}lienne } est un groupe alg{\'e}brique complet et connexe. On peut montrer qu'un tel groupe alg{\'e}brique est n{\'e}cessairement commutatif.
Un  {\it  groupe lin{\'e}aire} est un sous-groupe alg{\'e}brique de  $\GL_n$ pour un entier
positif 
$n$.

Un th{\'e}or{\`e}me de  Rosenlicht {\'e}tablit  que tout homomorphisme de groupes alg{\'e}briques d'une vari{\'e}t{\'e} ab{\'e}lienne dans un groupe
lin{\'e}aire ou d'un groupe lin{\'e}aire connexe dans une vari{\'e}t{\'e}
ab{\'e}lienne,
est constant.

Remarquons qu'un groupe alg{\'e}brique fini est lin{\'e}aire 
et que toute vari{\'e}t{\'e} ab{\'e}lienne admet des sous-groupes finis (de torsion). Donc l'hypoth{\`e}se
de connexit{\'e} est n{\'e}cessaire dans l'{\'e}nonc{\'e} ci-dessus.

Voici un th{\'e}or{\`e}me de  structure d\^u {\`a} Chevalley

\begin{theoreme}
Soit  $G$ un groupe alg{\'e}brique connexe. Il
existe un sous-groupe alg{\'e}brique normal, connexe  et lin{\'e}aire $L$ de $G$ tel que le quotient 
$A=G/L$ soit une vari{\'e}t{\'e} ab{\'e}lienne.
Ce  $L$ est unique et il contient tous les sous-groupes alg{\'e}briques lin{\'e}aires connexes
de $G$. 
\end{theoreme}

Les homomorphismes de groupes alg{\'e}briques surjectifs {\`a} noyaux finis sont appel{\'e}s
isog{\'e}nies. Pour les groupes alg{\'e}briques commutatifs 
sur un  corps  fini on a  un th{\'e}or{\`e}me
de structure plus fort. Si  $G$ est un groupe alg{\'e}brique connexe
commutatif sur un
corps fini $\FF_q$ et si  
 $1\rightarrow L\rightarrow G\rightarrow A\rightarrow 1$ est la suite exacte stricte donn{\'e}e
par le th{\'e}or{\`e}me de  Chevalley, alors il existe une isog{\'e}nie,
d{\'e}finie sur  $\FF_q$, de $G$ vers le produit direct $L\times A$.

\`A ma connaissance, tous les groupes alg{\'e}briques
utilis{\'e}s {\`a} ce jour en cryptographie sont, de fa{\c c}on plus ou moins visible,  des
vari{\'e}t{\'e}s jacobiennes g{\'e}n{\'e}ralis{\'e}es.

Un  {\it tore} sur le corps  $K$ est un groupe alg{\'e}brique connexe
commutatif $\TT$ de dimension $d$, 
qui devient isomorphe {\`a}
$(\Gm )^d$ apr{\`e}s
extension des scalaires de  $K$ {\`a} une extension s{\'e}parable  $L$. Un tel corps
$L$ est appel{\'e} corps de d{\'e}composition du tore $\TT$.
Donc un tore de dimension $d$ est un tordu du groupe alg{\'e}brique  $(\Gm )^d$
et  ces tores sont classifi{\'e}s par 
$H^1(K,\Aut_{K^s}((\Gm )^d))$. 

Observons que nous voulons tordre  $(\Gm )^d$ en tant que groupe
alg{\'e}brique et non seulement en tant que vari{\'e}t{\'e}.
Donc le groupe d'automorphismes $\Aut_{K^s}((\Gm )^d)$ 
qui nous pr{\'e}occupe est le groupe
des automorphismes de groupe alg{\'e}brique et non pas le groupe complet
des automorphismes de la vari{\'e}t{\'e}  $(\Gm )^d$.

Un  automorphisme $\agot$ du groupe alg{\'e}briques  $(\Gm )^d$ 
est d{\'e}crit par

$$\agot (g_1,g_2,...,g_d)=\agot ( (g_i)_{1\le i\le d} )=(\prod_{1\le j\le d}
g_j^{\alpha_{i,j}})_{1\le i\le d}.$$

Ici les  $\alpha_{i,j}$ sont des entiers tels que la  matrice 
$A=(\alpha_{i,j})$ ait un d{\'e}terminant {\'e}gal {\`a} $\pm 1$.

Ainsi le groupe $\Aut_{K^s}((\Gm )^d)$ est isomorphe {\`a}  $\GL_d(\ZZ)$ et l'action de 
Galois sur ce groupe est triviale. Donc  
$Z^1(K,\Aut_{K^s}((\Gm )^d))$ n'est autre que 
 $Hom(\GG_K,\Aut_K((\Gm )^d))$ et 
$H^1(K,\Aut_{K^s}((\Gm )^d))$ est le quotient de ce dernier
ensemble par $\Aut_K((\Gm )^d)$ agissant par conjugaison.

On suppose d{\'e}sormais que le  corps de base $K$ est fini
de caract{\'e}ristique $p$. Le groupe de Galois absolu
$\GG_K$ est procyclique et un {\'e}l{\'e}ment de $H^1(K,\Aut_{K^s}((\Gm )^d))$
est donn{\'e} par la classe de conjugaison 
de l'image de l'automorphisme  de Frobenius $F$.

Pour  $L\supset K$ une extension de degr{\'e}  $d$ de corps finis, on
identifie $(\Gm )^d$ au produit $\Pi_{\Gal(L/K)} {\Gm }$  de $d$ facteurs $\Gm $
indic{\'e}s
par les $K$-automorphismes de $L$.  

On note $\fgot_d$ l'automorphisme de $\Pi_{\Gal(L/K)} {\Gm }$ qui permute 
les composantes comme $F$ agit sur les 
indices dans 
 $\Gal(L/K)$ : 

$$\fgot_d(g_1,g_F,g_{F^2},..., g_{F^{d-1}})=(g_{F^{d-1}}, g_1,g_F,g_{F^2}, ...,
g_{F^{d-2}}).$$

On note $\cG_d$ le tordu de $\Pi_{\Gal(L/K)} {\Gm }$ associ{\'e} au cocycle
$F\mapsto \fgot_d$. C'est la restriction de Weil de $\Gm $ le long de $L/K$.
C'est un tore de  dimension $d$, d{\'e}fini sur $K$,
et  qui se d{\'e}compose sur  $L$.  Une de ses propri{\'e}t{\'e}s int{\'e}ressantes est
que $\cG_d(K)$ est  isomorphe,  en tant que  groupe,  {\`a} $\Gm (L)=L^*$.
Notons que le groupe alg{\'e}brique $\cG_d$ d{\'e}pend du corps
de base $K$ et du degr{\'e} $d$.

En effet, il existe un  $L$-isomorphisme

$$I :
(\Gm )^d=\prod_{\Gal(L/K)}\Gm \rightarrow \cG_d$$
\noindent  tel que 
${}^FI=I\circ \fgot_d$.

Pour tout diviseur  $a$ de $d$ posons  $d=ab$. Soit  $M$ 
l'extension de degr{\'e}  $a$
de  $K$. Alors $K\stackrel{a}{\subset}  M \stackrel{b}{\subset}
L$.  
Le groupe de Galois  $\Gal(L/K)$ 
est engendr{\'e} par le  Frobenius $F$ et le groupe de Galois
 $\Gal(L/M)$ est engendr{\'e}  par  $F^a$.

La restriction de $L$ {\`a} $M$ d{\'e}finit un {\'e}pimorphisme de groupes  
$\Gal(L/K)\rightarrow
\Gal(M/K)$. On en d{\'e}duit  l'existence
d'un homomorphisme  de groupes alg{\'e}briques $\nu_{a,d}$ de  $\cG_d$  dans
$\cG_a$.

En effet, on a un morphisme 

$$\cN_{a,d} :
(\Gm )^d=\prod_{\Gal(L/K)}\Gm \rightarrow
(\Gm )^a=\prod_{\Gal(M/K)}\Gm $$
\noindent  d{\'e}fini par

\begin{eqnarray*}
\cN_{a,d}(g_1,g_F,g_{F^2}, ... ,
g_{F^{d-1}})
&=&(g_1g_{F^a}g_{F^{2a}}...g_{F^{(b-1)a}},
g_Fg_{F^{a+1}}g_{F^{2a+1}}...g_{F^{(b-1)a+1}},\\
g_{F^2}g_{F^{a+2}}g_{F^{2a+2}}...g_{F^{(b-1)a+2}}, &...& , 
g_{F^{a-1}}g_{F^{a+a-1}}g_{F^{2a+a-1}}...g_{F^{(b-1)a+a-1}})
\end{eqnarray*}
\noindent et tel que $\cN_{a,d} \circ \fgot_d =\fgot_a \circ \cN_{a,d}$.

On note  $\cT_d$ l'intersection des noyaux des morphismes 
$\nu_{a,d}$ pour tous les diviseurs stricts $a$ de $d$.
C'est   un tore
de  dimension $\phi(d)$ tel que  $\cT_d(K)\subset
\cG_d(K)=\Gm (L)=L^*$ est le sous groupe des points qui ont  norme {\'e}gale {\`a} $1$ dans toute sous-extension 
de  $L/K$.
Ce sous groupe a pour cardinalit{\'e}  $\Phi_d(q)$ 
o{\`u}  $\Phi_d(X)$ est le  $d$-i{\`e}me
polyn{\^o}me cyclotomique et $q$ la cardinalit{\'e} de  $K$.

\subsubsection{Le tore de Lucas}

Soit encore $K$ un corps fini de caract{\'e}ristique impaire $p$.  Soit
$D\in K^*$ un scalaire qui n'est pas un carr{\'e} dans $K$. Soit $L=K(\sqrt D)$.
Soit $\AA^2/K$ le  plan  affine et 
 $\cU$ l'ouvert   d{\'e}fini par l'in{\'e}galit{\'e}  
$x^2-Dy^2\not =0$. On construit un
$L$-isomorphisme
de $\cU\otimes_KL$  dans  $(G_m)^2$ en envoyant
  $(x,y)$ sur $(x+y\sqrt
D , x-y\sqrt D)$. L'isomorphisme inverse est  $I : (G_m)^2\rightarrow
\cU\otimes_KL$ d{\'e}fini par 
$I(z_1,z_2)=(x,y)=(\frac{z_1+z_2}{2},\frac{z_1-z_2}{2\sqrt D})$.

On v{\'e}rifie  que  $I^{-1}{}^\sigma I$ est l'identit{\'e} si
$\sigma$ fixe $\sqrt D$ et l'application d'inversion sinon. Cela prouve
que 
$\cU$ est $K$-isomorphe {\`a} $\cG_2$.

Le sous groupe   $\cT_2\subset \cG_2$ est d{\'e}fini par la condition suppl{\'e}mentaire
que 
$z_1z_2=1$ ou de fa{\c c}on {\'e}quivalente  $x^2-Dy^2=1$.

Donc  $\cT_2$ est le ferm{\'e} de  $\AA^2$ d{\'e}fini par l'{\'e}quation
 $x^2-Dy^2=1$. Le tore $\cT_2$ est appel{\'e} tore de Lucas.

La loi de groupe sur  $\cG_2$ et  $\cT_2$ est donn{\'e}e par les applications de multiplication
et d'inversion :

\begin{equation*}\xymatrix{
\mu: & \cG_2\times \cG_2   \ar@{->}[r]&\cG_2 \\
 & ((x_1,y_1),(x_2,y_2))   \ar@{|->}[r]&(x_1x_2+Dy_1y_2, y_1x_2+x_1y_2)}
\end{equation*}

\noindent et 

\begin{equation*}\xymatrix{
i: & \cG_2  \ar@{->}[r]&\cG_2 \\
 & (x,y)   \ar@{|->}[r]&(x, -y)}
\end{equation*}

Le groupe $\cT_2(K)$ des points $K$-rationnels du tore de Lucas
est parfois pr{\'e}f{\'e}r{\'e} au groupe $\Gm(K)$ parce que le logarithme
discret y est suppos{\'e} un peu plus difficile. Notons
que les points de $\Gm(K)$ sont repr{\'e}sent{\'e}s par une
seule coordonn{\'e}e affine alors que  les points
de $\cT_2(K)$ sont repr{\'e}sent{\'e}s ici par leurs deux coordonn{\'e}es $x$ et $y$.
On a donc une  repr{\'e}sentation
  deux fois plus longue  pour $\cT_2(K)$ alors que les deux groupes
sont de tailles comparables. C'est un inconv{\'e}nient s{\'e}rieux.

Une premi{\`e}re parade, qui est assez g{\'e}n{\'e}rique,
consiste {\`a}  noter que l'exponentiation par un entier $k$ dans le groupe
alg{\'e}brique   $\cT_2$ est donn{\'e}e par $[k] : \cT_2\rightarrow \cT_2$ 
avec 

$$[k](x,y)=(\sum_{0\le 2l\le k} (x^2-1)^lx^{k-2l}\left(\begin{array}{c}
k\\2l \end{array}\right), y\sum_{0\le 2l+1\le k} (x^2-1)^lx^{k-2l}\left(\begin{array}{c}
k\\2l+1 \end{array}\right)).$$

En particulier, la coordonn{\'e}e $x$ de  $[k]P$ ne d{\'e}pend que de la coordonn{\'e}e 
$x$ de  $P$. Ceci simplement par ce que 
$[k]$ commute {\`a} l'inversion $i$.

Il est donc naturel de consid{\'e}rer la vari{\'e}t{\'e} quotient
$\cX_2=\cT_2/\{ 1,i  \}$. Celle-ci n'est autre
que la droite affine
avec $x$ pour coordonn{\'e}e. 

Ce quotient  n'est plus un groupe alg{\'e}brique mais il conserve
une action
du monoide multiplicatif des entiers positifs $(\ZZ_{>0},\times)$, donn{\'e}e par
les applications  $[k]$ d'exponentiation.

Cela suffit pour faire de la cryptographie {\`a} base de logarithme discret
pourvu que l'on se contente de l'exponentiation et que l'on renonce
{\`a} la multiplication.

Cette vari{\'e}t{\'e} $\cX_2$ pr{\'e}sente l'avantage 
d'{\^e}tre rationnelle  : les points sont d{\'e}crits par une seule
coordonn{\'e}e. On peut ainsi repr{\'e}senter des
probl{\`e}\-mes de logarithme discret dans le sous-groupe  de cardinal
$q+1$ de $L^*$ avec seulement $\log_2(q)$ bits (ceux qui suffisent {\`a} d{\'e}crire la coordonn{\'e}e
 $x$). 

C'est l'id{\'e}e {\`a} l'origine du syst{\`e}me LUC de  \cite{LUC}.

Il est de la premi{\`e}re importance pour cette
m{\'e}thode que la vari{\'e}t{\'e} $\cX_2$ soit  {\it rationnelle}.

Notons que le m{\^e}me proc{\'e}d{\'e} est
utilis{\'e} pour les courbes elliptiques car le quotient d'une courbe
elliptique par son involution est rationnel lui aussi.

Mais il y a mieux.  Rubin et Silverberg    rappellent  dans
\cite{rubinsilv, rubinsilv2} que 
le tore $\cT_2$ lui-m{\^e}me est  rationnel comme  $K$-vari{\'e}t{\'e}.

En effet, l'{\'e}quation  $x^2-Dy^2=1$ est rendue homog{\`e}ne en posant  $x^2-Dy^2=t^2$ qui est aussi
$x^2-t^2=Dy^2$ ou encore  $(x-t)(x+t)=Dy^2$.  On pose  $u=\frac{x-t}{y}$ et il vient
une param{\'e}trisation  $\frac{x}{y}=\frac{u+\frac{D}{u}}{2}$ et $\frac{t}{y}=\frac{-u+\frac{D}{u}}{2}$ donc une
 param{\'e}trisation de  $x^2-Dy^2=t^2$ par

\begin{eqnarray*}
x&=&u^2+D\\
t&=&-u^2+D\\
y&=&2u
\end{eqnarray*}
\noindent ce qui en coordonn{\'e}es affines donne $x=\frac{D+u^2}{D-u^2}$ et $y=\frac{2u}{D-u^2}$.

Ce qui se produit ici est que  $G_m$ est une sous vari{\'e}t{\'e} de  
$\PP^1$. Bien qu'il existe
un tordu non-trivial  du groupe  
$G_m$, le $1$-cocycle associ{\'e} s'annule dans  $H^1(K,\Aut_{K^s}(\PP^1))$.

Une cons{\'e}quence int{\'e}ressante est qu'un {\'e}l{\'e}ment de  $\cT_2$ peut {\^e}tre repr{\'e}sent{\'e} par une seule
coordonn{\'e}e $u$
et que la loi de groupe peut s'exprimer en terme de cette unique coordonn{\'e}e.

En effet, soit  $P_1$ de  $u$-coordonn{\'e}e $u_1$ et 
$P_2$ de $u$-coordonn{\'e}e  $u_2$, un calcul sans  myst{\`e}re
donne la       $u$-coordonn{\'e}e  de  $P_3$ que l'on note $u_3$ :

$$u_3= \frac{D(u_1+u_2)}{u_1u_2+D}.$$

L'{\'e}l{\'e}ment identit{\'e} de  $\cG_2$ a une   $u$-coordonn{\'e}e  {\'e}gale {\`a} $0$ et l'application 
d'inversion change   $u$ en $-u$.

Le point de coordonn{\'e}es  $x=-1$ et  $y=0$ correspond {\`a}  $u=\infty$.

Si  $u$ parcourt   $K\cup \{\infty\}$  il  repr{\'e}sente
les  $q+1$ points de $\cT_2(K)$.

 Silverberg (qui s'appelle Alice) et Rubin (qui ne s'appelle pas Bob)
 montrent que le cryptosyst{\`e}me XTR \cite{XTR} utilise 
une vari{\'e}t{\'e}
de type $\cX$ tout comme LUC. La question qui se pose alors est de 
d{\'e}terminer
dans quels cas cette vari{\'e}t{\'e} quotient d'un tore est rationnelle. 
Et si elle l'est,
de donner une param{\'e}trisation explicite. La m{\^e}me question se pose
pour les tores $\cT_d$ eux m{\^e}mes et elle est assez ouverte.
Le tore $\cT_d$ et son quotient $\cX_d$ sont  de dimension $\phi(d)$.
Pour $d=2\times3\times 5$ on   a $\phi(d)=8$ et  on cherche des 
param{\'e}trisations ...

Le  tore de Lucas est utilis{\'e} (au moins implicitement) depuis longtemps.
Il est le ressort de la m{\'e}thode appel{\'e}e ``$p+1$'' pour factoriser
des entiers naturels en produit de facteurs premiers.

D'une mani{\`e}re g{\'e}n{\'e}rale, factoriser un nombre entier $N$ revient
{\`a} calculer le nombre de points $\ZZ/\NN\ZZ$ rationnels d'un groupe
alg{\'e}brique $\bG$ bien choisi. C'est {\'e}vident si on choisit
$\bG=\Gm$ mais il peut {\^e}tre plus habile de choisir $\bG=\cT_2$
le tore de Lucas ou encore $\bG=A$ une vari{\'e}t{\'e} ab{\'e}lienne.

Dans le m{\^e}me ordre d'id{\'e}e, on prouve qu'un nombre $P$  est premier
en prouvant que $\bG(\ZZ/P\ZZ)$ a le cardinal  attendu o{\`u} $\bG$
est un groupe alg{\'e}brique bien choisi. 
Si $\bG=\Gm$ on attend $P-1$, si $\bG=\cT_2$
on attend $P+1$, et si $\bG$ est une vari{\'e}t{\'e} ab{\'e}lienne {\`a} multiplication
complexe, on sait aussi {\`a} quoi s'attendre, gr{\^a}ce {\`a} la th{\'e}orie
de Shimura.

La vari{\'e}t{\'e} (pas si grande) des groupes alg{\'e}briques commutatifs
a donc {\'e}t{\'e} explor{\'e}e largement par de nombreux auteurs int{\'e}ress{\'e}s
{\`a} l'une ou l'autre de ces questions : factorisation, primalit{\'e}, 
logarithme discret
protocoles cryptographiques.

On trouve dans  \cite{Cohen} un {\'e}tat de  ces m{\'e}thodes  et dans 
\cite{Serre} une introduction aux groupes alg{\'e}briques.

\subsection{Des exemples d'instances faibles  du logarithme
  discret}\label{subsection:araki}

Le probl{\`e}me du logarithme discret n'est pas toujours difficile.
On a vu qu'il est facile dans le groupe additif d'un corps fini.
Il est facile aussi dans un groupe $C=\prod_iC_i$ produit direct de petits
groupes. Dans ce cas, il se d{\'e}compose en autant de probl{\`e}mes 
de logarithmes discrets dans les facteurs $C_i$. On {\'e}vite donc
les groupes $\Gm(\Fq)$ si $q-1$ est friable (produit de petits facteurs
premiers).

Dans ce paragraphe  nous montrons un autre exemple d'instance
faible du logarithme discret. La cause,  ici,   est analytique :
on r{\'e}duit le logarithme discret {\`a} un logarithme $p$-adique.

L'alin{\'e}a  \ref{subsubsection:Riesel} pr{\'e}sente la m{\'e}thode
de Riesel pour calculer le logarithme discret dans le groupe
multiplicatif $\Gm(\ZZ/p^k\ZZ)$.  L'alin{\'e}a
\ref{subsubsection:sas} expose l'extension de cette
 m{\'e}thode au cas des courbes elliptiques de trace $1$, d'apr{\`e}s
les travaux de Smart, Araki, Satoh, Semaev.

\subsubsection{Rappels sur les logarithmes  $p$-adiques}

Un {\'e}l{\'e}ment de l'anneau  $\Zp$ des entiers
 $p$-adiques est not{\'e}

$$z=z_0+z_1p+z_2p^2+...+z_{k-1}p^{k-1}+O(p^k)$$ 
\noindent o{\`u} les $z_i$ sont des entiers 
tels que $0\le z_i<p$ et $k$ est la pr{\'e}cision absolue requise.

On note $v_p$ la valuation $p$-adique sur $\Qp$.
On sait que $\Qps=<p>\times \Zps$. On note $\U=\Zps$
et pour tout entier $n\ge 1$ soit     $\U_n=1+p^n\Zp$.
La r{\'e}duction modulo $p$ donne une suite exacte

\begin{equation}\label{equation:exact}
1\rightarrow \U_1 \rightarrow \Zps {\rightarrow }\FF^*_p
\rightarrow 1.
\end{equation}

Le lemme de  Hensel montre que $\U$ contient
le groupe  $\V$  des racines  $(p-1)$-i{\`e}mes de l'unit{\'e}. Donc la suite
exacte ci-dessus se d{\'e}compose  et on a $\U=\V \times  \U_1.$

Reste {\`a} d{\'e}crire $\U _1$.
Si  $p\ge 3$ on construit un homomorphisme
de groupes  topologiques entre 
 $(\U _1,\times)$ et   $(p\Zp,+)$.

Si  $p=2$  alors  $\U _1=<-1>\times \U _2$ et on construit
un homomorphisme
de groupes  topologiques entre 
  $(\U _2,\times)$ et   $(4\Zp,+)$.

Dans les deux cas, l'isomorphisme est donn{\'e} par la s{\'e}rie
logarithme 

$$\Log (z)= -\sum_{n=1}^{\infty}{\frac{(-1)^n}{n}}(z-1)^n.$$

L'application inverse est l'exponentielle

$$\exp (z)=\sum_{n=0}^\infty {\frac{z^n}{ n!}}.$$

La s{\'e}rie   $\Log (z)$ converge sur 
$\U_1$ ($\U_2$ si $p=2$) et la s{\'e}rie 
$\exp (z)$  converge sur   $p\Zp$ ($4\ZZ _2$ si
$p=2$). Donc la structure de   $\Qps$ est 

$$\ZZ\times \ZZ/(p-1)\ZZ \times \Zp$$
\noindent pour   $p$ impair et

$$\ZZ \times \ZZ/2\ZZ \times \ZZ _2$$
\noindent pour $p=2$.

\subsubsection{Calcul de logarithmes discrets dans $(\ZZ/p^k\ZZ)^*$}\label{subsubsection:Riesel}

Nous supposons dans la suite que $p$ est un entier premier impair.

Soit  $k\ge 2$ un entier et
 $b_k$ un  g{\'e}n{\'e}rateur de $(\ZZ/p^k\ZZ)^*$. Soit
 $c_k$ un autre {\'e}l{\'e}ment de  $(\ZZ/p^k\ZZ)^*$. On veut calculer
le logarithme discret $\ell_k=\log_{b_k}c_k$. On
voit  $(\ZZ/p^k\ZZ)^*$ comme le quotient de 
$\Zps$ par  $\U_k$ et on fixe 
un rel{\`e}vement  $b$  de  $b_k$ dans
$\Zps$
et un rel{\`e}vement   $c$  de  $c_k$ dans
$\Zps$.

Soient  $c_1$ et  $b_1$ les images de  $c$ et $b$ dans 
$(\ZZ/p\ZZ)^*=\Zps/\U_1$ et supposons que le logarithme
$\ell_1=\log_{b_1}c_1\in \ZZ/(p-1)\ZZ$ est connu.
On identifie  la classe de congruence 
$\ell_1$ modulo $p-1$ {\`a} son repr{\'e}sentant dans l'intervalle  $[0,p-1[$.

Cet entier  peut {\^e}tre calcul{\'e} par recherche
exhaustive au prix de $O(p)$ op{\'e}rations dans  $\ZZ/p\ZZ$.

On pose  $C=cb^{-\ell_1}$ et on v{\'e}rifie que  $C\in \U_1$. Soit
 $B=b^{p-1}\in \U_1-\U_2$.

On calcule {\`a} l'aide du d{\'e}veloppement en s{\'e}rie convergente 

$$L=\frac{\Log C}{\Log B} \pmod{p^{k-1}}$$
\noindent
donc  $C=B^L\pmod{p^k}$ et  $c=b^{\ell_1}b^{(p-1)L}\pmod{p^k}$
\noindent ce qui nous donne le logarithme discret cherch{\'e}

$$\log_{b_k}c_k=\ell_1+(p-1)L\pmod{p^{k-1}(p-1)}.$$

Ainsi, le calcul du logarithme discret
dans  $\ZZ/p^k\ZZ$ se r{\'e}duit au calcul d'un logarithme discret dans 
 $\ZZ/p\ZZ$. Pour  $p$ fix{\'e}  et $k$ tendant vers l'infini, on obtient
un exemple de grand groupe multiplicatif o{\`u} le logarithme discret 
se calcule en temps polynomial en le logarithme de la taille $S=(p-1)p^{k-1}$ du groupe.

\subsubsection{Un exemple}

Soit   $p=3$ et   $k=10$. On choisit

$$b = 59045 =2 + 3 + 2.3^2  + 2.3^3  + 2.3^4  + 2.3^5  + 2.3^6  +
2.3^7  + 2.3^8  + 2.3^9  + O(3^{10})$$
\noindent un g{\'e}n{\'e}rateur de   $\ZZ/3^{10}\ZZ$.

Posons   $B=b^2= 1 + 2.3 + 3^2  + O(3^{10})$. C'est un g{\'e}n{\'e}rateur
de   $\U_1$.

Soit  

$$c=24731= 2 + 2.3 + 2.3^2  + 2.3^4  + 2.3^5  + 2.3^7  + 3^9  +
O(3^{10}).$$

On veut calculer  $\log _{b_{10}} c_{10}$.

Puisque   $c = b \pmod 3$ on a   $\ell _1=1$ et on pose

$$C=c/b =1 + 2.3 + 3^3  + 2.3^4  + 3^6  + 2.3^7  + 2.3^8  + 3^9
+ O(3^{10}).$$

On calcule   $L=\Log (C)/\Log (B)$ gr{\^a}ce au d{\'e}veloppement 
en s{\'e}rie du logarithme 

$$L= 1 + 3 + 2.3^2  + 2.3^3  + 3^5  + 3^7  + O(3^9)$$

Donc   $L=2506 \pmod {3^9}$ et   $\ell _{10}= 1+2\times 2506=5013$.

Dans l'article original \cite{Riesel}de Riesel, le logarithme
$p$-adique est remplac{\'e} par 
le quotient de  Fermat.

\subsubsection{Courbes elliptiques sur un corps local}

On suppose encore que $p$ est premier impair.
Soit  $\cE$  une courbe elliptique d'{\'e}quation affine

$$E : y^2+a_1xy+a_3y=x^3+a_2x^2+a_4x+a_6$$
\noindent avec 
$a_1, a_2, a_3, a_4, a_6$ dans $\Zp$.

On suppose que  $\cE$ a bonne r{\'e}duction
: le discriminant $\Delta$ est une unit{\'e}. On note $\tilde \cE$
la r{\'e}duction de $\cE$ modulo $p$.

La r{\'e}duction modulo $p$ d{\'e}finit un {\'e}pimorphisme de groupes
$\rho : \cE(\Qp)\rightarrow
\tilde \cE(\Fp)$
de $\cE(\Qp)$ vers  $\tilde \cE(\Fp)$. 
Le noyau de $\rho$ est form{\'e} des points proches $p$-adiquement
de l'origine. On note   $z=-x/y$ le param{\`e}tre local
en l'origine et pour tout $k$ positif on note 
 $\cE_k(\Qp)$ l'ensemble des points $P\in \cE(\Qp)$ tels que 
$v_p(z_P)\ge k$. On a la suite exacte

\begin{equation}
0\rightarrow \cE_1(\Qp)\rightarrow \cE(\Qp)\stackrel{\rho {}}{\rightarrow}
\tilde \cE(\Fp)\rightarrow 0.
\end{equation}
\noindent qui est tr{\`e}s proche de la suite \ref{equation:exact}.

On {\'e}tudie donc la structure de groupe de 
$\cE_1(\Qp)$. La loi de groupe sur $\cE_1(\Qp)$
est exprim{\'e}e en terme du param{\`e}tre $z$. En effet $z(P+Q)$
est une s{\'e}rie formelle en  $z(P)$ et  $z(Q)$. Cette s{\'e}rie formelle n'est autre
que le 
groupe formel associ{\'e} {\`a} $\cE$, c'est-{\`a}-dire le d{\'e}veloppement de Taylor {\`a}
l'origine
de la loi d'addition. 
Les coefficients de ce d{\'e}veloppement sont des polyn{\^o}mes en les coefficients
de l'{\'e}quation de $\cE$ : 

\begin{equation*}
z(P+Q)=z(P)+z(Q)
-a_1z(P)z(Q)-a_2(z(P)^2z(Q)+z(Q)^2z(P))+\cdots
\end{equation*}

Le param{\`e}tre local $z$ induit une bijection 
de  $\cE_1(\Qp)$ sur  $p\ZZ_p$. La s{\'e}rie $F(z_1,z_2)$ converge sur
 $p\ZZ_p\times
p\ZZ_p$ et pour tout  $P$ et  $Q$ dans   $\cE_1(\Qp)$ on a $z(P+Q)=F(z(P),z(Q))$.
On pose  $z_1 \pf z_2=F(z_1,z_2)$ pour  $z_1$ et  $z_2$ dans  $p\ZZ_p$. Cela
fait
de  $z$ un homomorphisme de groupes topologiques 
de  $(\cE_1,+)$ vers $(p\Zp,\oplus_\cF)$.

Il existe un    {\it  logarithme formel }  associ{\'e} au groupe formel $F$ et
not{\'e} $\Logf (z)$. Il est caract{\'e}ris{\'e} par l'identit{\'e}

$$\Logf (F(z_1,z_2))=\Logf (z_1)+ \Logf (z_2).$$

La s{\'e}rie r{\'e}ciproque de  $\Logf (z)$ est not{\'e}e   $\expf (z)$.

Un simple calcul montre que $\Logf (z)=\sum _{n=1}^\infty {\frac{b_n}{ n}}z^n$
et $\expf (z)=\sum _{n=1}^\infty {\frac{c_n}{n!}}z^n$ o{\`u} les $b_n$ et   $c_n$
sont dans  $\Zp$, et  $b_1=c_1=1$.
On en d{\'e}duit que  $\Logf (z)$ converge pour 
$v_p(z)>0$ et   $\expf (z)$ pour   $v_p(z)>1/(p-1)$.

En consid{\'e}rant l'isomorphisme compos{\'e} 

$$\cE_1(\Qp) \stackrel{z}{\rightarrow} (p\ZZ_p,\oplus_\cF)
\stackrel{\Logf}{\rightarrow}(p\Zp,+)$$
\noindent on  prouve que    $\cE_1(\Qp)$ est sans torsion et que la r{\'e}duction
modulo $p$  est injective sur le sous-groupe de torsion de $\cE(\Qp)$. On note  $\tilde o= \#\tilde \cE (\Fp)$.
Si   $\tilde o$ est premier {\`a}  $p$  alors   $\tilde \cE (\Fp)$ se rel{\`e}ve en un
groupe
de torsion dans   $\cE (\Qp)$ et la r{\'e}duction modulo  $p$ induit  une
bijection 
entre les torsions de 
$\cE(\Qp)$  et   $\tilde \cE (\Fp)$. Donc

$$\cE(\Qp)\sim \cE_1(\Qp)\times \tilde \cE (\Fp)\sim \Zp \times \tilde \cE 
(\Fp).$$

Supposons maintenant que   $\tilde o=\cE(\Fp)=p$. C'est un cas tr{\`e}s
particulier car la trace de l'endomorphisme de Frobenius
de $\cE$ vaut $1$. Alors 
$\cE(\Qp)$ est coinc{\'e} dans la suite exacte

$$0\rightarrow \Zp \rightarrow \cE(\Qp)\rightarrow \ZZ/p\ZZ\rightarrow
0.$$

Donc il est isomorphe {\`a}   $\Zp$ (cas I) ou a  $\Zp\times
\ZZ/p\ZZ$ (cas II).

Supposons que l'on se trouve dans le cas I.
Donc  $\cE(\Qp)$ est sans torsion et il existe une bijection
entre  $\tilde \cE (\Fp )$ et
 $\cE_1/\cE_2$.
Soit en effet   $\tilde P$ un   point dans   $\tilde \cE (\Fp)$ et soit 
 $P\in \cE (\Qp)$ un rel{\`e}vement   $\tilde P$. Le   point $[p]P$ est dans  $\cE_1 (\Qp)$ car
$\widetilde {[p]P}=[p]\tilde P=0$. En outre, si nous choisissons un autre
rel{\`e}vement  $P'$ de $\tilde P$, alors 
$[p]P-[p]P'=[p](P-P')\in \cE _2$. On a donc une application

$$\Pi : \tilde \cE (\Fp) \rightarrow \cE_1/\cE_2.$$

Elle est injective. En effet, soit   $P\in \cE(\Qp)$ tel que 
$[p]P\in \cE_2$. La multiplication par
 $p$ d\'efinit une   bijection entre   $\cE_1$ et $\cE_2$. Il y a un   $R\in
 \cE_1$ tel que   $[p]R=p[P]$.
Donc $[p](R-P)=0$ et puisque on est dans le cas  I on a   $P=R\in \cE_1$ et  $\tilde P=0$.

\subsubsection{La m{\'e}thode de  Smart-Araki-Satoh-Semaev}\label{subsubsection:sas}

Pour calculer le logarithme discret dans 
une courbe elliptique de trace $1$, ils
 utilisent la bijection  $\Pi$ de l'alin{\'e}a pr{\'e}c{\'e}dent
et   transforment  un probl{\`e}me de  logarithme
discret dans   $\tilde \cE (\Fp)$ en un probl{\`e}me de  logarithme dans  $\cE_1(\Qp)$. 
Ce dernier logarithme est un logarithme elliptique et se calcule efficacement
en raison de ses propri{\'e}t{\'e}s analytiques.

Il reste {\`a} s'assurer que l'on se trouve dans le cas I de l'alin{\'e}a pr{\'e}c{\'e}dent. En fait Voloch  a not{\'e}
que
le cas  II corresponds au cas o{\`u}  $\cE$ est le rel{\`e}vement canonique
de $\tilde \cE$ modulo $p^2$. Comme ce rel{\`e}vement canonique est unique, il
n'est pas
tr{\`e}s difficile {\`a} {\'e}viter\ldots

D'un point de vue pratique, le d{\'e}veloppement du logarithme elliptique
est obtenu en int{\'e}grant la forme diff{\'e}rentielle canonique $\omega$.

On d{\'e}veloppe   $x$, $y$ et   $\omega $
en  $z=-x/y$ et on calcule 

\begin{eqnarray*}
x&=&z^{-2}-a_1z^{-1}-a_2 .... \\
y&=& -x/z\\
\omega &=& {\frac{dy}{3x^2+2a_2x+a_4-a_1y}}=(1+a_1z+(a_1^2+a_2)z^2+...)dz\\
\Logf &=&\int \omega =z+{\frac{a_1}{ 2}}z^2+{\frac{a_1^2+a_2}{ 3}}+...
\end{eqnarray*}

Finissons par un exemple. 
Soit   $p=655387895585476301924777$ et 
$\tilde \cE$ la courbe d'{\'e}quation

$$\tilde E : y^2+xy=x^3+ 114287067913850793676921x+
349073807889941681395769$$
\noindent sur  $\Fp$.

Soit   $\tilde P=(170219448,14643735815400225272219)$ 
un  point dans $\tilde \cE(\Fp)$. On v{\'e}rifie que 
$[p]\tilde P=0$.

Soit  $\tilde Q=(71434243993450257115004,316317604915944437378529)$ un autre
point de   $\tilde \cE(\Fp)$.

On choisit un relev{\'e}  $\cE$ de   $\tilde\cE$
sur  $\Qp$. Par exemple on choisit  $\cE$ d'{\'e}quation 

$$ E : y^2+xy=x^3+ 114287067913850793676921x+
349073807889941681395769$$
\noindent sur  $\Qp$.

On cherche un  $P=(x_P,y_P)$ dans   $\cE(\Qp)$
au dessus de   $\tilde P$. On fixe   $x_P=170219448$ par exemple et on r{\'e}sout
dans   $\Qp$ 
l'{\'e}quation en   $y$ 

$$y^2+x_Py=x_P^3+114287067913850793676921x_P+
349073807889941681395769.$$

On choisit la racine congrue {\`a}

$14643735815400225272219$ modulo $p$.

Donc 

$$y_P=14643735815400225272219 + 241062303587335366096866.p +O(p^2).$$

Alors 

\begin{eqnarray*}
\Pi (\tilde P)&=&[p]P\\
&=&\scriptstyle (246304660834813598643589.p^{-2}  + O(p^{-1}),
213491610344127745612815.p^{-3}  + O(p^{-2})).\\
\end{eqnarray*}

 Et 

$$\Logf ([p]P)=z_P+O(z_P^2)=
-\frac{246304660834813598643589}{213491610344127745612815}p+O(p^2).$$

De m{\^e}me

$$\Logf ([p]Q)=169836480309236709243708.p+O(p^2).$$

Le quotient est 
$\log _P (Q)=123456789\pmod p$.


On voit que les courbes elliptiques qui ont $p$ points rationnels sur $\Fp$
sont
des instances  faibles pour le logarithme discret.
Cette observation a {\'e}t{\'e} faite
par Araki  et Satoh \cite{Araki},  Smart \cite{Smart}, Semaev 
\cite{Semaev} ind{\'e}pendemment. 
Ruck \cite{Ruck} et Voloch ont unifi{\'e} et g{\'e}n{\'e}ralis{\'e}  ces travaux dans le langage naturel
de la cohomologie galoisienne.

\section{Isog{\'e}nies et cryptographie}\label{section:isog}

Dans cette section, je  d{\'e}cris un EHD qui n'est
pas un probl{\`e}me de logarithme discret. Lorsque j'ai pr{\'e}sent{\'e} cet exemple
en 1997 au s{\'e}minaire de cryptographie de l'ENS, il s'agissait d'une curiosit{\'e}
(voir  \cite{couv}).
Mais les progr{\`e}s r{\'e}alis{\'e}s dans le calcul explicite des isog{\'e}nies et les
travaux men{\'e}s par  Charles, Lauter, Jao et Venkatesan depuis lors,
montrent qu'il n'est pas irr{\'e}aliste de fonder la s{\'e}curit{\'e} d'un cryptosyst{\`e}me
sur la difficult{\'e} de trouver un morphisme entre deux objets.

Le paragraphe \ref{subsection:myh2s} d{\'e}crit l'action du groupe des classes
 d'un ordre
quadratique
sur les courbes elliptiques {\`a} multiplication complexe par cet ordre, d'un point
de vue algorithmique. Cette situation produit un candidat EHD.

Le paragraphe \ref{subsection:graph} rappelle quelques  propri{\'e}t{\'e}s des  graphes
d'isog{\'e}nies  entre courbes elliptiques et donne une id{\'e}e de leur int{\'e}r{\^e}t
cryptographique. 

\subsection{L'espace homog{\`e}ne des courbes ordinaires {\`a} multiplication
par un ordre quadratique}\label{subsection:myh2s}

Soit  $\Fq$ un corps de cardinal $q=p^d$ et $E$ une courbe elliptique
sur  $\Fq$, suppos{\'e}e ordinaire. L'anneau des endomorphismes de $E$
est isomorphe {\`a} un ordre quadratique  $\OO$. On fixe un tel isomorphisme
$\iota : \End(E)\rightarrow \OO$.
On note  $t$ la trace de l'endomorphisme de  Frobenius $\Phi$. Donc $\# E(\Fq)=q+1-t.$

On pose  $\Delta =t^2-4q$  et on suppose que $\Delta$ est sans facteur
carr{\'e}. Donc $\OO=\ZZ[\Phi]$ est maximal. 

Si  $\agot $ est un id{\'e}al de  $\OO$ premier {\`a} $p$,
on note $\Ker \agot \subset E$
l'intersection des noyaux des isog{\'e}nies de $E$ appartenant {\`a} $\iota^{-1}(\agot)$.

Soit  $F$ le quotient de $E$ par $\Ker \agot$ et 
$I_\agot: E\rightarrow F$  l'isog{\'e}nie quotient. Soit $\kappa : \End(F)\rightarrow
\OO$ l'isomorphisme d{\'e}fini par
 $\kappa(\alpha)= \iota (I_\agot^{-1} \alpha I_\agot)$.

On note $[\agot]$ la classe de $\agot$ dans $\Pic(\cO)$ et
on pose $[\agot] . (E,\iota)= (F,\kappa)$. On d{\'e}finit ainsi 
une action du groupe des
classes
$G=\clk$ de $\OO$ sur l'ensemble  des classes d'isomorphismes de couples 
$(E,\iota)$ form{\'e}s d'une courbe elliptique sur $\Fq$ et d'un isomorphisme de
$\End(E)$ vers $\OO$.

Cette action admet deux orbites, permut{\'e}es par la conjugaison complexe. Soit
$H$ l'une des deux orbites. L'action de $G$ sur $H$ est simplement transitive. Le cardinal $S=\# G = \# H$ est le
nombre de classes de $\OO$. On a $S=\Delta^{\frac{1}{2}+o(1)}$.

Il est facile de v{\'e}rifier qu'une courbe elliptique 
$E$ a un anneau d'endomorphismes isomorphe {\`a} $\OO$. Il suffit de v{\'e}rifier que
 $\# E =p+1-t$ {\`a} l'aide de la m{\'e}thode de Schoof \cite{Schoof} ou plus vite encore si on dispose d'une factorisation
de $p+1-t$ et si ce dernier entier est sans facteur carr{\'e}.

Si  $\ell\not = p$ est un entier premier impair d{\'e}compos{\'e} dans $\OO$, le polyn{\^o}me $f(X)=X^2-tX+q$
a deux racines distinctes   $\lambda$ et  $\mu$ modulo $\ell$  et on a 
$\ell =\lgot \bar\lgot $ avec $\lgot  =(\ell, \Phi -\lambda)$ et
$\bar\lgot  =(\ell, \Phi -\mu )$.

Le noyau $\Ker \lgot$ est un sous-groupe de $E[\ell]$. Il correspond {\`a} un
facteur de degr{\'e} $\frac{\ell-1}{2}$ du polyn{\^o}me de $\ell$-division. Ce facteur
se calcule en temps polyn{\^o}mial en $\log q$ et $\ell$ par divers moyens (soit
brutalement, soit en suivant les id{\'e}es de Schoof, Atkin, Elkies et quelques
autres \cite{Schoof2}).

Le quotient de $E$ par $\Ker \lgot$ se calcule en temps polyn{\^o}mial en 
$\ell$ 
et $\log q$ {\`a} l'aide des formules de V{\'e}lu
ou, plus efficacement encore, en utilisant des m{\'e}thodes plus r{\'e}centes
propos{\'e}es par Elkies, Schoof, Lercier, moi-m{\^e}me, Morain, Salvy,
Schost, etc.
On peut consulter \cite{schost} qui est un texte r{\'e}cent sur le sujet.
Notons que le calcul d'isog{\'e}nies se d{\'e}compose en deux {\'e}tapes : trouver
d'abord le noyau, puis quotienter la courbe, ou bien au contraire, trouver 
d'abord la courbe
quotient 
({\`a} l'aide d'{\'e}quations modulaires) puis  en d{\'e}duire le noyau.

Un nombre premier sur deux se d{\'e}compose dans $\OO$.
Si cette proportion est respect{\'e}e 
pour les petits nombres premiers,  en admettant l'hypoth{\`e}se de Riemann
g{\'e}n{\'e}ralis{\'e}e, on dispose
d'un ensemble d'{\'e}l{\'e}ments de $G$ dont   l'action  sur $H$ se calcule en temps
polyn{\^o}mial en $\log q$ et qui engendrent $G$.

Une combinaison de ces {\'e}l{\'e}ments {\`a} coefficients entiers al{\'e}atoires et assez
grands  produit un {\'e}l{\'e}ment al{\'e}atoire de $G$ avec distribution assez proche de
la distribution uniforme (sur la composante connexe du point de d{\'e}part).

Il est raisonnable de penser  qu'il n'existe pas d'algorithme rapide pour
calculer l'unique  classe de $\clk$ qui envoie une courbe $E$ sur une courbe  $F$.
On a donc un candidat EHD s{\'e}rieux qui ne provient pas du logarithme discret.

On ne conna{\^\i}t pas d'algorithme  pour calculer $S$ le nombre de classes
de $\OO$ en temps polyn{\^o}mial en $\log \Delta$. Mais conna{\^\i}tre le cardinal
exact d'un EHD n'est pas indispensable aux protocoles que nous avons pr{\'e}sent{\'e}s.

Pour construire un EHD tel que d{\'e}crit dans le paragraphe pr{\'e}c{\'e}dent,
on choisit d'abord un corps fini $\Fq$. On choisit une courbe elliptique
au hasard et on calcule la trace $t$ de l'endomorphisme de Frobenius avec
la m{\'e}thode de Schoof. On v{\'e}rifie que $\Delta=t^2-4q$ est
sans facteur carr{\'e} et se factorise ais{\'e}ment (sinon on recommence).

On collectionne alors les petits nombres premiers $\ell$ qui se d{\'e}composent
dans $\ZZ[\Phi]$.

\subsection{Graphes d'isog{\'e}nies}\label{subsection:graph}

Dans la situation du paragraphe pr{\'e}c{\'e}dent, on fixe un r{\'e}el
positif $\delta$ et on pose $B=(\log q)^{2+\delta}$. Soit $\cG$
le graphe dont les sommets sont les {\'e}l{\'e}ments de $H$ et dont les cot{\'e}s sont les
paires $\{(E,\iota),(F,\kappa)\}$ d'{\'e}l{\'e}ments de $H$ telles qu'il existe un id{\'e}al
premier 
$\agot$
de $\OO$, de degr{\'e} d'inertie  $1$,  de norme $\le B$ et tel que  $[\agot].E = F$. C'est un graphe de
Cayley. Il est $k$-r{\'e}gulier o{\`u}  $k$ est le nombre d'id{\'e}aux premier de degr{\'e}
d'inertie $1$ et de norme $\le B$ dans $\OO$.
On note $M=[m_{e,f}]_{e,f \in H}$ sa matrice d'adjacence. Donc $m_{e,f}=1$ si
$\{e,f\}$
est un cot{\'e} du graphe, et $m_{e,f}=0$ sinon. 

Le matrice $\frac{1}{k}M$ est une matrice de Markov, correspondant {\`a} une
marche al{\'e}atoire dans le graphe  : si l'on se trouve au   sommet $e$ du graphe,
on choisit un des $k$ cot{\'e}s  issus de $e$ avec probabilit{\'e} uniforme et on
avance  le long de ce sommet.

La distribution uniforme de probabilit{\'e}s sur l'ensemble $H$ des sommets est un
vecteur propre de $M$ et sa valeur propre est $k$.

Toutes les valeurs propres de $M$ sont r{\'e}elles et de valeur absolue $\le
k$. On peut  consulter le petit livre \cite{Sarnak} de Sarnak sur les graphes
et les formes modulaires pour toutes ces questions.

Les valeurs propres non-triviales (celles dont la valeur absolue est $<k$)
ralentissent la convergence du processus de Markov vers    la distribution
uniforme. Si elles sont petites, alors cette convergence est rapide.

Jao,  Miller et Venkatesan montrent que pour le graphe ci-dessus, sous reserve que
l'hypoth{\`e}se
de Riemann g{\'e}n{\'e}ralis{\'e}e soit correcte, le graphe est connexe et  les valeurs
propres
non-triviales ont une valeur absolue $O(k^{\beta})$ pour tout  $\beta >
\frac{1}{2} + \frac{1}{\delta +2}$. Comme $\delta $ est positif, ces valeurs
propres
sont nettement s{\'e}par{\'e}es de la valeur propre associ{\'e}e {\`a} la distribution
uniforme.
Donc le processus markovien converge vite.

On dit que la famille des graphes ainsi construits est une famille de  graphes
d'expansion.

Classiquement, on construit  plut{\^o}t des graphes d'isog{\'e}nies {\`a} l'aide de courbes
supersinguli{\`e}res, comme dans   \cite{pizer}.  En effet, les matrices d'adjacences, appel{\'e}es matrices de
Brandt, expriment l'action d'op{\'e}rateurs de Hecke sur des espaces de formes
modulaires.
Leurs valeurs propres sont  major{\'e}es {\`a} l'aide de la conjecture de
Ramanujan.  Ces graphes r{\'e}guliers, appel{\'e}s graphes de Pizer, sont encore meilleurs que
les pr{\'e}c{\'e}dents, car leurs valeurs propres non-triviales sont en valeur absolue $\le 2\sqrt{k-1}$. On dit que ce
sont des graphes de Ramanujan.

Outre leurs remarquables propri{\'e}t{\'e}s spectrales, les graphes d'isog{\'e}nies
pr{\'e}sentent un grand int{\'e}r{\^e}t calculatoire. Ce sont de grands graphes dans
lesquels
on peut circuler facilement (cela revient {\`a} calculer des isog{\'e}nies de petit
degr{\'e}). Mais il est difficile de trouver un chemin entre deux sommets donn{\'e}s,
comme nous l'avons vus dans le paragraphe \ref{subsection:myh2s}. Charles, Goren et
Lauter proposent dans \cite{lauter} d'utiliser  cette propri{\'e}t{\'e} pour construire des fonctions de
hachage
cryptographique. Le principe est le suivant : on 
appelle $\cC$ l'ensemble des entiers premiers
$\le B$ qui se d{\'e}composent dans $\OO$ et on forme un ensemble
$\cB$ en choisissant pour tout $\ell$ dans $\cC$ un id{\'e}al au dessus
de $\ell$. On choisit une origine parmi les sommets du graphe et on fixe
une bijection entre les  lettres de l'alphabet et les  id{\'e}aux
  dans  $\cB$,  de sorte qu'{\`a} tout mot de longueur quelconque on peut 
associer
un chemin dans le graphe. Le sommet o{\`u} l'on aboutit {\`a} l'issue de ce
cheminement est la valeur de la fonction de hachage. Trouver deux mots
qui se hachent sur le m{\^e}me sommet revient {\`a} trouver un cycle non-trivial
dans le graphe.

\bibliography{couveignes-smf}

\end{document}